\documentclass[onefignum,onetabnum]{siamonline250211}   

\usepackage{graphicx} 
\usepackage[utf8]{inputenc}
\usepackage{cite} 
\usepackage{amssymb} 
\usepackage{amsmath}
\usepackage{hyperref}
\usepackage{enumitem}

\usepackage{makecell} 

\renewcommand{\fbox}{} 

\long\def\beginpgfgraphicnamed#1#2\endpgfgraphicnamed{\includegraphics{#1}}

\newcommand{\R}{\mathbb{R}}

\title{Visibility of heteroclinic networks}

\author{Sofia B. S. D. Castro,\thanks{Faculdade de Economia and Centro de Matem\'atica, Universidade do Porto, Portugal
  (\email{sdcastro@fep.up.pt})}
\and Claire M. Postlethwaite\thanks{Department of Mathematics, University of Auckland, New Zealand
  (\email{c.postlethwaite@auckland.ac.nz})}
\and Alastair M. Rucklidge\thanks{School of Mathematics, University of Leeds, UK
  (\email{A.M.Rucklidge@leeds.ac.uk})}}

\begin{document}

\maketitle

\begin{abstract}
The concept of stability has a long history in the field of dynamical systems: stable invariant objects are the ones that would be expected to be observed in experiments and numerical simulations. Heteroclinic networks are invariant objects in dynamical systems associated with intermittent cycling and switching behaviour, found in a range of applications. In this article, we note that the usual notions of stability, even those developed specifically for heteroclinic networks, do not provide all the information needed to determine the long-term behaviour of trajectories near heteroclinic networks. To complement the notion of stability, we introduce the concept of \emph{visibility}, which pinpoints precisely the invariant objects that will be observed once transients have decayed. We illustrate our definitions with examples of heteroclinic networks from the literature. 
\end{abstract}

\begin{keywords}
Stability; Heteroclinic networks; Heteroclinic cycles.
\end{keywords}

\begin{MSCcodes}
34C37, 34D05, 34D20, 34D45, 37C81
\end{MSCcodes}


\section{Introduction}\label{sec:intro}

The concept of \emph{stability} is an old one in dynamical systems and is invariably associated with the persistence and/or observability of an object or outcome. 
The work of Poincar\'e~\cite{Poincare1885} and of Lyapunov~\cite{Liapounoff1907} are unavoidable landmarks, but many other authors have made important contributions. 
We refer the interested reader to the historical accounts provided by Leine~\cite{Leine2010}, Roque~\cite{Roque2011}, Mawhin~\cite{Mawhin1994} and references therein. 
In the context of dynamical systems, we typically think of an invariant object as being \emph{stable} when either nearby solutions do not move far away (Lyapunov stability~\cite{Liapounoff1907}) or nearby solutions converge to the object (Milnor attractor~\cite{Milnor1985}). 
When the two occur simultaneously, we talk about \emph{asymptotic stability} (see Auslander \textit{et al.}~\cite{Auslander1964}).

In this paper we introduce the concept of \emph{visibility}, to complement the classic notion of stability. We illustrate this with examples that have arisen in the study of heteroclinic networks.
These are invariant objects that are finite unions of heteroclinic cycles, themselves the union of finitely many equilibria and trajectories connecting them in a cyclic way.

A systematic study of the asymptotic stability of heteroclinic cycles dates back to the work of Krupa and Melbourne~\cite{Krupa1995}. 
A heteroclinic cycle which is a proper subset of a heteroclinic network can never be asymptotically stable even though it can be highly attracting. 
This is due to the fact that near an equilibrium that belongs to more than one cycle there are always some trajectories that move away from at least one of these cycles. In particular, initial conditions which lie exactly on a connecting orbit to an equilibrium \emph{not} in the cycle of interest will have this property.
The notions of \emph{essential asymptotic stability} (e.a.s.), and of \emph{fragmentary asymptotic stability} (f.a.s.), were introduced by Melbourne~\cite{Melbourne1991} and by Podvigina~\cite{Podvigina2012} respectively, to address this. 
A heteoclinic cycle or network is f.a.s.\ if it satisfies the conditions for asymptotic stability for points in a set of positive measure contained in a neighbourhood of the cycle or network, while it is e.a.s.\ if this set of points tends to a set of full measure as the size of the neighbourhood tends to zero. This is in contrast to the requirements of asymptotic stability, which are to be satisfied for \emph{all} points in a neighbourhood of the invariant set of interest.

Hence, ordered by decreasing order of `stability', heteroclinic cycles or networks can be asymptotically stable, essentially asymptotically stable, or fragmentarily asymptotically stable. 
All these definitions also require the conditions for Lyapunov stability, namely that trajectories do not move too far away. 
If this fails, the invariant set can be   \emph{quasi-asymptotically stable}~\cite{Glendinning1994} or a \emph{Milnor attractor}~\cite{Milnor1985}. 
In both cases, trajectories are attracted to the invariant set as time tends to infinity although they may move far away during a transient.

The motivation for this work appears when relating the stability of an invariant set to its observability: there are `very stable' heteroclinic networks that can never appear entirely in numerical simulations or as the outcome of an experiment. 
On the other hand, there are heteroclinic cycles classified as `not very stable' that invariably appear in numerical simulations. 
A quintessential example is the Kirk--Silber (KS) network~\cite{Kirk1994}, discussed in detail in Subsection~\ref{sec:KS} below. 
It is possible to choose parameters such that the entire network is e.a.s.\ but only one sub-cycle of the network is visible in simulations. For a different set of parameters, in which one cycle is e.a.s.\ but the other is not even f.a.s., it is this latter cycle that is visible (in the long term) in simulations. 
This illustrates the need for an alternative notion which provides precise information about what we can expect to observe in simulations and experiments. 
We argue that what we define as \emph{visibility}, and its various nuances, precisely pinpoints the invariant objects we can expect to observe.

This gap between a somewhat stable invariant set and its observability is not a strictly theoretical issue as we illustrate with examples from game theory and from population dynamics. 
That stability occasionally provides ambiguous information is clear in the population dynamics of the Rock--Paper--Scissors--Spock--Lizard (RPSSL) game studied by the authors in~\cite{Postlethwaite2022} and~\cite{Castro2022a} as well as those of the jungle game with four species studied in~\cite{Castro2024}. 
In both games, species compete amongst themselves in a non-transitive manner (also known as \emph{cyclic dominance}) and the dynamics correspond to those near a heteroclinic network. 
When only some of the species survive, a heteroclinic cycle involving the equilibria corresponding to these surviving species is expected as the outcome of simulations. 
However, under the current definitions of stability, the whole network of the RPSSL game can be  asymptotically stable even when as few as three species survive. 
More complex outcomes not captured by the stability of the network are found in~\cite{Postlethwaite2022}.
In the jungle game, although the whole network is asymptotically stable there is always one species that becomes extinct.
These two games in population dynamics illustrate the ambiguity contained in stating that the whole network is asymptotically stable in terms of survival of the species involved.

Therefore, the main purpose of this paper is to establish a concept that fills a gap left by the current definitions of stability. Visibility is more related to outcomes than stability, especially in the context of heteroclinic dynamics. The tools to establish visibility are the same we use for stability.

In the next section we present and discuss some of the notions of stability in the literature, and in Section~\ref{sec:examples} we can illustrate how these fail to provide complete information about what is observed in numerical experiments. 
These examples are known but our presentation includes new numerical simulations and presents previously unobserved/unreported behaviour in some cases.
We use the examples to motivate the variations we introduce in the concept of visibility.
Visibility is defined in Section~\ref{sec:visibility}. An invariant set that is visible always appears in numerical simulations in its entirety.
For instance, the KS network briefly mentioned above is never visible. 
In the population dynamics examples, only the cycles corresponding to the surviving species are visible, while the network is not.
We provide a discussion and illustration of visibility that clearly distinguishes it from any of the available notions of stability. We end with a discussion in the last section.

\section{Preliminaries}

In this section we give background on commonly used definitions of stability, and also on  heteroclinic cycles and networks. 
Let $X$ be a compact flow-invariant set for the dynamics described by the ODE 
\begin{equation}\label{eq:ODE}
\dot{x} = f(x)
\end{equation}
for $x \in \R^n$ and $f: \R^n \rightarrow \R^n$ smooth.   Recall that an \emph{equilibrium solution} $\xi_0$ satisfies $f(\xi_0)=0$. 

\begin{definition}\label{def:hetero-cycle}
	A \emph{heteroclinic cycle} is a finite collection of equilibria $\{\xi_1, \dots, \xi_m\}$ of~\eqref{eq:ODE}, together with a set of heteroclinic connecting trajectories $\{\gamma_1(t),\dots, \gamma_m(t)\}$, where $\gamma_j(t)$ is a solution of~\eqref{eq:ODE} such that  $\lim_{t\rightarrow -\infty} \gamma_j(t)=\xi_j$, $\lim_{t\rightarrow \infty} \gamma_j(t)=\xi_{j+1}$ and $\xi_{m+1}\equiv\xi_1$.
\end{definition}

A \emph{heteroclinic network} is a connected union of finitely many heteroclinic cycles. 
The literature contains numerous different definitions of heteroclinic cycles and networks: subtleties occur because there may be equilibria in the network with two-dimensional unstable manifolds, and whether or not to include the entire unstable manifold in the heteroclinic network can result in different definitions. 
In this paper, we use simply the definition given at the start of this paragraph, and for each example, we are specific about the object under investigation. Given Definition~\ref{def:hetero-cycle}, a heteroclinic cycle is a closed flow-invariant set.

More generally, heteroclinic cycles and networks may connect invariant objects more complicated than equilibria, such as periodic orbits~\cite{Aguiar2006,Kirk2008} or chaotic sets~\cite{Ashwin1997f}, but we do not consider these possibilities here. 
In generic systems, heteroclinic cycles and networks are of high co-dimension, but when~\eqref{eq:ODE} contains invariant subspaces, then they may exist for open sets of parameters values, that is, they are \emph{robust}.
An early example of a robust heteroclinic cycle can be found in~\cite{Guckenheimer1988} and is addressed in Subsection~\ref{sec:GH} below.

Several notions of stability have been developed to deal with the properties of compact flow-invariant sets, in particular with heteroclinic cycles and networks. We present next those that are most common.

Let $\varphi(t,x)$ denote the flow of equation~\eqref{eq:ODE}: the point in the orbit of $x$ at time $t$.
The \emph{$\omega$-limit} set of a point $x \in \R^n$ is 
$$
\omega(x) = \{ y \in \R^n: \mbox{there is an increasing sequence of times }t_n \rightarrow +\infty \mbox{ with }\!\! \lim_{n \rightarrow +\infty}\! \varphi(t_n,x) = y  \}.
$$

Our first two definitions encompass two important properties of stability: \emph{Lyapunov stability} requires that trajectories starting near an invariant set~$X$ stay near $X$, and \emph{quasi-asymptotic stability} requires that trajectories reach~$X$ in the limit as $t$ increases to infinity.

\begin{definition}[adapted from Kuznetsov~\cite{Kuznetsov1998}]\label{def:Lyapunov-stability}
    We say that an invariant set~$X$ is \emph{Lyapunov stable} if for any neighbourhood $U$ of $X$ there exists a neighbourhood $V$ of $X$ such that for all $x \in V$ and $t>0$,
    $$
    \varphi(t,x) \in U.
    $$
\end{definition}

If $d(a,b)$ denotes the distance between two points $a$ and $b$, then the distance between a point $a$ and a set $B$ is defined as
$$
d(a,B) = \inf \{ d(a,b): \; b \in B \}.
$$

\begin{definition}[adapted from Glendinning~\cite{Glendinning1994}]\label{def:quasi-stability}
We say an invariant set~$X$ is \emph{quasi-asymptotic stable} iff there exists $\delta>0$ such that if $d(y,X)<\delta$ then $d(\varphi(y,t),X)\rightarrow 0$ as $t\rightarrow \infty$.
\end{definition}

Similar to quasi-asymptotic stability, is the concept of a \emph{Milnor attractor}. 
\begin{definition}[Milnor~\cite{Milnor1985}]\label{def:attractor-Milnor}
    A closed subset $X \subset \R^n$ will be called an \emph{attractor} if it satisfies two conditions:
    \begin{enumerate}
        \item the \emph{realm of attraction}, $\rho(X)$, consisting of all points $x \in \R^n$ for which $\omega(x) \subseteq X$, must have strictly positive measure; and
        \item there is no strictly smaller closed set $X^\prime \subset X$ so that $\rho(X^\prime)$ coincides with $\rho(X)$ up to a set of measure zero.
    \end{enumerate}
\end{definition}

\emph{Asymptotic stability} is usually defined to be the combination of Lyapunov and quasi-asymptotic stability~\cite{Glendinning1994,Kuznetsov1998}.
There are different definitions but they coincide when the flow-invariant set $X$ is closed.
See~\cite{Podvigina2020} for a discussion.
 
Let $V_\delta(X)$ represent the set of points within the state space at a distance at most $\delta$ from~$X$.
Podvigina~\cite{Podvigina2012} defines \emph{$\delta$-basin of attraction} of a compact flow-invariant set $X$ as
$$
    \mathcal{B}_{\delta} (X) = \{x \in \R^n: \; \lim_{t \rightarrow +\infty} d(\varphi(t,x), X)=0  \mbox{  and  } \varphi(t,x) \in V_{\delta}(X)\}.
$$
The \emph{basin of attraction} of $X$ is just
$\mathcal{B} (X) = \{x \in \R^n: \; \lim_{t \rightarrow +\infty} d(\varphi(t,x), X)=0\}$.
This is the same as the realm of attraction in Milnor's definition (Definition~\ref{def:attractor-Milnor}).

We observe that the basin of attraction is not necessarily a neighbourhood of~$X$, and in the cases of heteroclinic cycles and networks, the basins frequently have thick or thin cusps removed~\cite{Kirk1994}.
Removing cusps leads to the definitions of essential asymptotic stability (Definition~\ref{def:eas}, which removes only thin cusps from the neighbourhood), and fragmentary asymptotic stability (Definition~\ref{def:fas}, which removes any type of cusp).
Essential asymptotic stability is the strongest notion of stability after asymptotic stability.
It originates in Melbourne~\cite{Melbourne1991} but is usually stated in the words of Brannath as follows.

\begin{definition}[Definition 1.1 in Brannath~\cite{Brannath1994}]\label{def:as-relative}
    Given any subset $N$ of $\R^n$, a closed invariant subset $X \subset \bar{N}$ ($\bar{N}$ denotes the closure of $N$) is said to be \emph{asymptotically stable, relative to the set~$N$} if for every neighbourhood $U$ of $X$ there is a neighbourhood $V$ of $X$, such that
    $$
    \forall \; x \in V \cap N , \; \varphi(t,x) \in U \mbox{  and   } \omega(x) \subseteq X.
    $$
\end{definition}
Here, asymptotically stable relatively to the set~$N$ means that initial conditions restricted to be in~$N$ will stay close to~$X$ and will asymptote to~$X$, but $N$ does not need to be a neighbourhood of~$X$.
The next definition, essential asymptotic stability, requires that the relative measure of~$N$ goes to~$1$ as $X$ is approached.
\begin{definition}[Definition 1.2 in Brannath~\cite{Brannath1994}] \label{def:eas}
A closed invariant set $X$ is called \emph{essentially asymptotically stable} (\hbox{e.a.s.}) if it is asymptotically stable relative to a set $N \subset \R^n$ with the property that
\begin{align*}
\lim\limits_{\delta \to 0} \frac{\ell(V_{\delta}(X) \cap N)}{\ell(V_\delta(X))} = 1,
\end{align*}
where $\ell$ is the Lebesgue measure.
\end{definition}

\begin{definition}[Podvigina~\cite{Podvigina2012}]\label{def:fas}
    We call a compact invariant set $X$ \emph{fragmentarily asymptotically stable} (\hbox{f.a.s.}) if, for any $\delta > 0$,  $\ell(\mathcal{B}_{\delta} (X))>0$.
\end{definition}

Note that the definitions of f.a.s.\ and e.a.s.\ demand that trajectories starting near $X$ do not go far away before converging to~$X$.
Without this restriction, f.a.s.\ would be the same as part~1 in Definition~\ref{def:attractor-Milnor} (Milnor attractor).

\section{Illustrative heteroclinic cycles and networks}\label{sec:examples}

In this section we present a few examples that illustrate how the existing notions of stability fail to provide complete information about what to expect when running simulations or experiments.
We work in $\R^n$ for $n \geq 3$ since this is the lowest dimension that supports the existence of a robust heteroclinic cycle.
The simplest way of constructing heteroclinic cycles is to have reflection symmetry in all the coordinates, so that the equations are equivariant under $x_j\rightarrow-x_j$, for $j=1,\dots,n$, and all the coordinate hyperplanes are invariant.
Heteroclinic cycles are found in systems with more elaborate symmetries as well, for example, mode interactions with spherical symmetry~\cite{Armbruster1991}, but we consider only examples with reflection symmetry here.

The examples in this section all have equations of the form
\begin{equation} \label{eq:odegen}
	\dot{x}_j = x_j\left(1- \chi + \sum_{k=1}^n a_{kj} x_k^2\right),\quad j=1,\dots, n,
\end{equation}
where $\chi=\sum_{j=1}^n{x_j^2}$, and the $a_{kj}$ are parameters with $a_{kk}=0$.
This system has equilibria on the axes with a single coordinate equal to $\pm 1$, and all other coordinates equal to zero.
We label the equilibria on the $x_j$ axis as $\pm\xi_j$. Each parameter $a_{kj}$ is the eigenvalue of the Jacobian matrix at $\pm\xi_j$ in the $x_k$ direction.
The remaining eigenvalue (in the direction $x_j$) is equal to $-2$, and does not affect the dynamics, due to the Invariant Sphere Theorem~\cite{Field1989}.
The reflection symmetry of~\eqref{eq:odegen} ensures that each coordinate axis and coordinate (hyper-) plane is invariant under the flow, and thus the positive orthant (which we write as $\R_+^n$) is also invariant under the flow.

For~\eqref{eq:odegen} restricted to $\R_+^n$, one can make the coordinate transformation $y_j= x_j^2$, and a rescaling of time, to get the equations
\begin{equation} \label{eq:odegenlv}
	\dot{y}_j = y_j\left(1- \chi + \sum_{k=1}^n a_{kj} y_k\right),\quad j=1,\dots, n.
\end{equation}
where $\chi=\sum_{j=1}^n{y_j}$.
These equations are the \emph{Lotka--Volterra} equations~\cite{Hofbauer1998} and are used in population dynamics contexts, where negative values of the coordinates have non-physical meanings. These equations are not equivariant, but still contain the same flow invariant coordinate axes and (hyper-) planes.

\subsection{The Guckenheimer--Holmes cycle}\label{sec:GH}

The prototypical example of a robust heteroclinic cycle is the \emph{Guckenheimer--Holmes} cycle~\cite{Guckenheimer1988}, also known as the \emph{Rock--Paper--Scissors} cycle. The governing equations are three-dimensional, and often written as:
\begin{align}
	\dot{x}_1&= x_1(1-\chi-c_{21}x_2^2+e_{31}x_3^2), \nonumber \\
		\dot{x}_2&= x_2(1-\chi-c_{32}x_3^2+e_{12}x_1^2), \label{eq:RPS} \\
			\dot{x}_3&= x_3(1-\chi-c_{13}x_1^2+e_{23}x_2^2). \nonumber
\end{align}
where $\chi=\sum_{j=1}^3{x_j^2}$, and $c_{jk}>0$ and $e_{jk}>0$ are parameters. 

Guckenheimer and Holmes imposed a three-fold symmetry $(x_1,x_2,x_3)\rightarrow (x_2,x_3,x_1)$, meaning there were only two distinct parameters. 
We do not impose that restriction here.

\begin{figure}
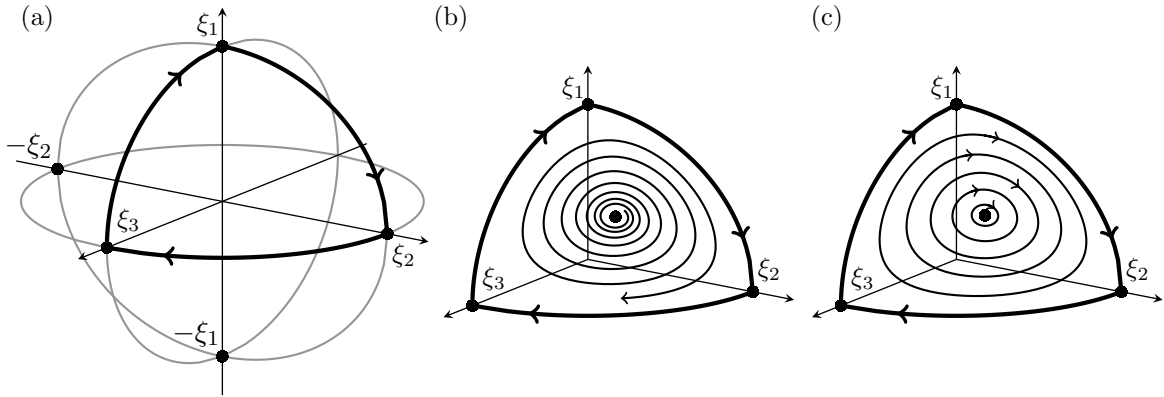

\hbox to \hsize{\hfil%
\fbox{\beginpgfgraphicnamed{GH_network_combined}%
\endpgfgraphicnamed}\hfil}

\caption{(a) Sketch of the Guckenheimer--Holmes heteroclinic `network', with the cycle in the positive orthant highlighted. (b) Sketch of the dynamics in the positive orthant when $\rho_{123}>1$. (c) Sketch of the dynamics in the positive orthant when $\rho_{123}=1$; note that there exists an infinite family of periodic orbits.
\label{fig:GHcycnet}
}
\end{figure}

The signs of the parameters ensure the existence of heteroclinic connections from $\pm\xi_j$ to $\pm\xi_{j+1}$ (indices taken modulo~3), lying within the (invariant) $(x_j,x_{j+1})$ coordinate plane.
Often the cycle in the positive orthant connecting the positive equilibria is considered to be the `Guckenheimer and Holmes cycle'.
It is also possible to consider the set of six equilibria and twelve connecting heteroclinic orbits as a heteroclinic network, formed from the union of eight heteroclinic cycles (see Figure~\ref{fig:GHcycnet}(a)). 
As mentioned above, the parameters $c_{jk}$ and $e_{jl}$ are related to the eigenvalues at the equilibria. 
Each equilibrium $\xi_j$ is a saddle point. 
The \emph{contracting} eigenvalue is equal to $-c_{jk}$ in the direction $x_k$. 
The \emph{expanding} eigenvalue is equal to $e_{jl}$ in the direction $x_l$.
In this way, the notation in \eqref{eq:odegen} relates to that in \eqref{eq:RPS} by $a_{kj}=e_{jk}$ and $a_{kj}=-c_{jk}$.
 
\begin{figure}
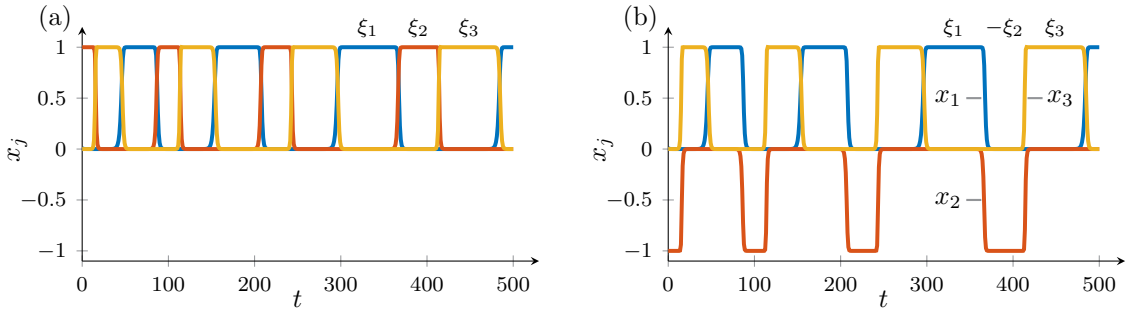

\hbox to \hsize{\hfil%
\fbox{\beginpgfgraphicnamed{x1x2x3_gh_timeseries_1}%
\endpgfgraphicnamed}\hfil%
\fbox{\beginpgfgraphicnamed{x1x2x3_gh_timeseries_2}%
\endpgfgraphicnamed}\hfil}     

\caption{Trajectories near the Guckenheimer--Holmes cycle~\eqref{eq:RPS}. In (a) the initial conditions have positive values for all coordinates, and the trajectory approaches the cycle between $+\xi_1$, $+\xi_2$ and $+\xi_3$. In (b), the initial condition has a negative value of $x_2$ and the trajectory asymptotes onto a cycle between $+\xi_1$, $-\xi_2$ and $+\xi_3$. Parameters are the same in both panels and are given in table~\ref{tab:parametervalues} in Appendix A.
\label{fig:GH_example}}
\end{figure}

It can be shown (by constructing an appropriate Poincar\'e map, see e.g.~\cite{Melbourne1989a}
for details), that the Guckenheimer and Holmes heteroclinic \emph{network} is asymptotically stable if
\[
\rho_{123}\equiv \frac{c_{13}c_{21}c_{32}}{e_{12}e_{23}e_{31}}>1.
\]
However, any trajectory will only converge to a single subcycle within the network. Each subcycle of this network cannot be asymptotically stable, since initial conditions near that subcycle but in the `wrong' orthant will move away. In Figure~\ref{fig:GH_example}, we show two trajectories near the Guckenheimer and Holmes cycle with different initial conditions. In panel (a), the initial condition lies in the positive orthant, and the trajectory asymptotes onto the cycle between $+\xi_1$, $+\xi_2$ and $+\xi_3$. In panel (b), the initial condition has a negative value of $x_2$ and the trajectory asymptotes onto a cycle between $+\xi_1$, $-\xi_2$ and $+\xi_3$. 
In fact, each subcycle is only f.a.s.\ if $\rho_{123}>1$. The basins of attraction of each subcycle are one of the eight orthants, including the boundary planes. Note that this means that the basins of attraction of each cycle have non-empty intersection.

Generally this discussion is avoided by either identifying the equilibria $\xi_j$ with $-\xi_j$ (invoking symmetry), or restricting the flow to the positive orthant~\cite{Krupa1997}. Then there is no longer a network of connections, but only a single cycle.
Most authors consider this cycle to be `the Guckenheimer and Holmes cycle', and that it is asymptotically stable when $\rho_{123}>1$.


However, there are some initial conditions for which the trajectories do not asymptote onto the entire cycle, namely those which lie within the invariant coordinate planes. 
These will instead converge to an equilibrium.
Hence, even though the cycle is classified as asymptotically stable some choices of initial conditions show only one of the equilibria in numerical simulations.
These initial conditions form only a set of measure zero, which may explain why this issue has not yet been addressed.

The equivalent Lotka--Volterra equations modelling intransitive competition between three species are 
\begin{align}
	\dot{y}_1&= y_1(1-\chi-c_{21}y_2+e_{31}y_3), \nonumber \\
		\dot{y}_2&= y_2(1-\chi-c_{32}y_3+e_{12}y_1), \label{eq:RPSLV} \\
			\dot{y}_3&= y_3(1-\chi-c_{13}y_1+e_{23}y_2). \nonumber
\end{align}
where $\chi=\sum_{j=1}^3{y_j}$. In these equations there is a single heteroclinic cycle, equivalent to that shown in Figure~\ref{fig:GHcycnet}(b). However, trajectories which start outside of $\R_+^3$ will diverge to $-\infty$; there is no network of connections as in the equivariant case. Thus, when the phase space is considered to be $\R^n$, the cycle cannot be asymptotically stable, but can be f.a.s., again if $\rho_{123}>1$. 

For the remainder of this paper, we will consider the dynamics of equations~\eqref{eq:odegen} or~\eqref{eq:odegenlv} with phase space as the positive orthant ($\R_+^n$). Then~\eqref{eq:odegen} and~\eqref{eq:odegenlv} are equivalent under a coordinate transformation and we make no further distinction between the Guckenheimer and Holmes network and the Guckenheimer and Holmes cycle.

Other subtleties regarding stability can occur when $\rho_{123}=1$. At this point we say that the Guckheimer and Holmes cycle is \emph{at resonance}. 
In system~\eqref{eq:RPS} then for these parameter values, it can be shown that an infinite family of periodic orbits exists (see Figure~\ref{fig:GHcycnet}(c))---there is a Hopf bifurcation from an equilibrium with three non-zero components at the same parameter values.
The Guckenheimer--Holmes cycle is Lyapunov stable; any initial condition close to the cycle will remain on a periodic orbit close to the cycle and thus remain a bounded distance from the cycle. This distance can be made smaller by choosing an initial condition closer to the cycle.

This example is one of several that motivate our definition of \emph{asymptotic visibility} in Section~\ref{sec:visibility}. 

\begin{figure}
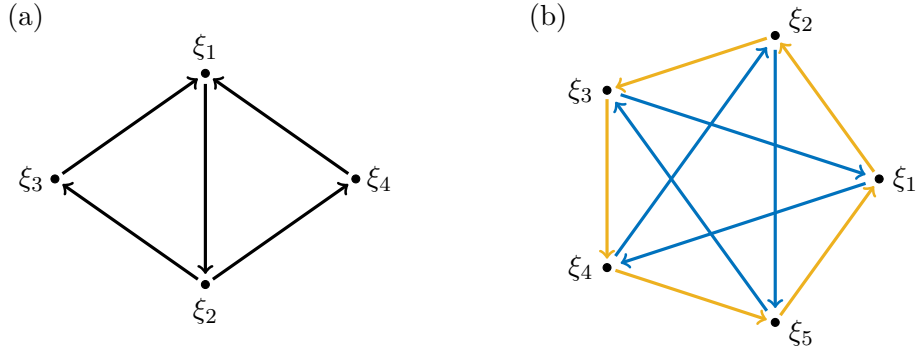

\hbox to \hsize{\hfil%
\fbox{\beginpgfgraphicnamed{KSnet_straight}%
\endpgfgraphicnamed}\hfil%
\fbox{\beginpgfgraphicnamed{RPSSL_pentacle_network}%
\endpgfgraphicnamed}\hfil}
\caption{(a)~The Kirk--Silber network, and 
(b)~the Rock--Paper--Scissors--Spock--Lizard network, both represented as directed graphs. Dots indicate equilibria, which lie on co-ordinate axes, and arrows represent heteroclinic connections. In (a) all connections are one-dimensional and contained in the same plane as the equilibria they connect.
The colours in (b) indicate two different types of connection: type A (amber) is a two-dimensional connection with one connecting trajectory in the plane containing the equilibria it connects and the remaining connecting trajectories in three-dimensional space; type B (blue) are one-dimensional connecting trajectories in the plane containing the equilibria.
See Section~\ref{sec:LS} and Figure~\ref{fig:2Dwu} for more detail on the type A (amber) connections.
\label{fig:KS_and_RPSSL_networks}}
\end{figure}

\subsection{The Kirk--Silber network}\label{sec:KS}

The Kirk--Silber network~\cite{Kirk1994} is formed by joining two Guckenheimer--Holmes cycles along an edge in $\R^4$ (see Figure~\ref{fig:KS_and_RPSSL_networks}(a)). Although equations were not given by Kirk and Silber, the network exists in the following system:
\begin{align}
	\dot{x}_1&= x_1(1-\chi-c_{21}x_2^2+e_{31}x_3^2+e_{41}x_4^2), \nonumber \\
		\dot{x}_2&= x_2(1-\chi+e_{12}x_1^2-c_{32}x_3^2-c_{42}x_4^2), \label{eq:KS} \\
			\dot{x}_3&= x_3(1-\chi-c_{13}x_1^2+e_{23}x_2^2-t_{43}x_4^2), \nonumber \\
	\dot{x}_4&= x_4(1-\chi-c_{14}x_1^2+e_{24}x_2^2-t_{34}x_3^2). \nonumber
\end{align}
where $\chi=\sum_{j=1}^4{x_j^2}$, and $c_{jk}>0$, $e_{jk}>0$ and $t_{jk}>0$ are parameters. 

As per the discussion in the previous section, if we consider the phase space to be $\R^4$, no cycle or network can be asymptotically stable because of the issues of trajectories in the `wrong' orthant. 
However, we will ignore those issues and restrict the phase space to $\R_+^4$. Even with this restriction, the Kirk--Silber system is heteroclinic network, not a cycle. Our main focus of attention is the relative stability of the cycles $\xi_1$-$\xi_2$-$\xi_3$ and $\xi_1$-$\xi_2$-$\xi_4$, and the network formed by their union. 

When considering the stability of the network or its subcycles, parameter combinations of interest are $\rho_{123}\equiv \frac{c_{13}c_{21}c_{32}}{e_{12}e_{23}e_{31}}$ and $\rho_{124}\equiv \frac{c_{14}c_{21}c_{42}}{e_{12}e_{24}e_{41}}$.
These denote, as usual, the ratio between the product of the absolute value of contracting and the product of the expanding eigenvalues at the nodes of a cycle.
These ratios determine whether trajectories approach the cycle, when the dynamics are restricted to the invariant subspace containing that cycle.
There are another six parameter combinations $\nu_{ijkl}$ which determine whether the cycles are stable to perturbations in the transverse direction. These are, for the $\xi_1$-$\xi_2$-$\xi_3$ cycle,
\begin{align*} 
\nu_{1234}&=\frac{c_{14}}{e_{12}}-\frac{c_{13}e_{24}}{e_{12}e_{23}}+\frac{c_{21}c_{13}t_{34}}{e_{12}e_{31}e_{23}},\\
\nu_{2314}&=-\frac{e_{24}}{e_{23}}+\frac{c_{21}t_{34}}{e_{23}e_{31}}+\frac{c_{14}c_{32}c_{21}}{e_{12}e_{23}e_{31}}, \\
\nu_{3124}&=\frac{t_{34}}{e_{31}}+\frac{c_{32}c_{14}}{e_{12}e_{31}}-\frac{c_{32}c_{13}e_{24}}{e_{23}e_{12}e_{31}},
\end{align*}
and the corresponding quantities for the $\xi_1$-$\xi_2$-$\xi_4$ cycle have the $3$ and $4$ switched in all indices. 
If any of $\nu_{ijk4}<0$, then the $\xi_1$-$\xi_2$-$\xi_3$ cycle is unstable to perturbations in the $x_4$ direction and similarly if any of $\nu_{ijk3}<0$, then the $\xi_1$-$\xi_2$-$\xi_4$ cycle is unstable to perturbations in the $x_3$ direction.

Kirk and Silber~\cite{Kirk1994} show that if $\rho_{123}>1$ and all $\nu_{ijk4}>0$ then the $\xi_1$-$\xi_2$-$\xi_3$ cycle is f.a.s.\ (although they did not use that term).
Similarly, if $\rho_{124}>1$ and all the $\nu_{ijk3}>0$, then the $\xi_1$-$\xi_2$-$\xi_4$ cycle is f.a.s. 
These sets of conditions can be satisfied simultaneously, and then the Kirk--Silber network is f.a.s. 
Kirk and Silber also show that the \emph{network} can be f.a.s., even if one of the cycles is not f.a.s. 
For instance, if one of the $\nu_{ijk3}<0$ but all $\nu_{ijk4}>0$, then trajectories which start close to the $\xi_1$-$\xi_2$-$\xi_3$ cycle will eventually `switch' to the $\xi_1$-$\xi_2$-$\xi_4$ cycle. 
However, they also show that it is not possible to choose parameters where trajectories switch in both directions.
Therefore, there are no initial conditions whose trajectories repeatedly visit the entire network, for any set of parameter values.

\begin{figure}
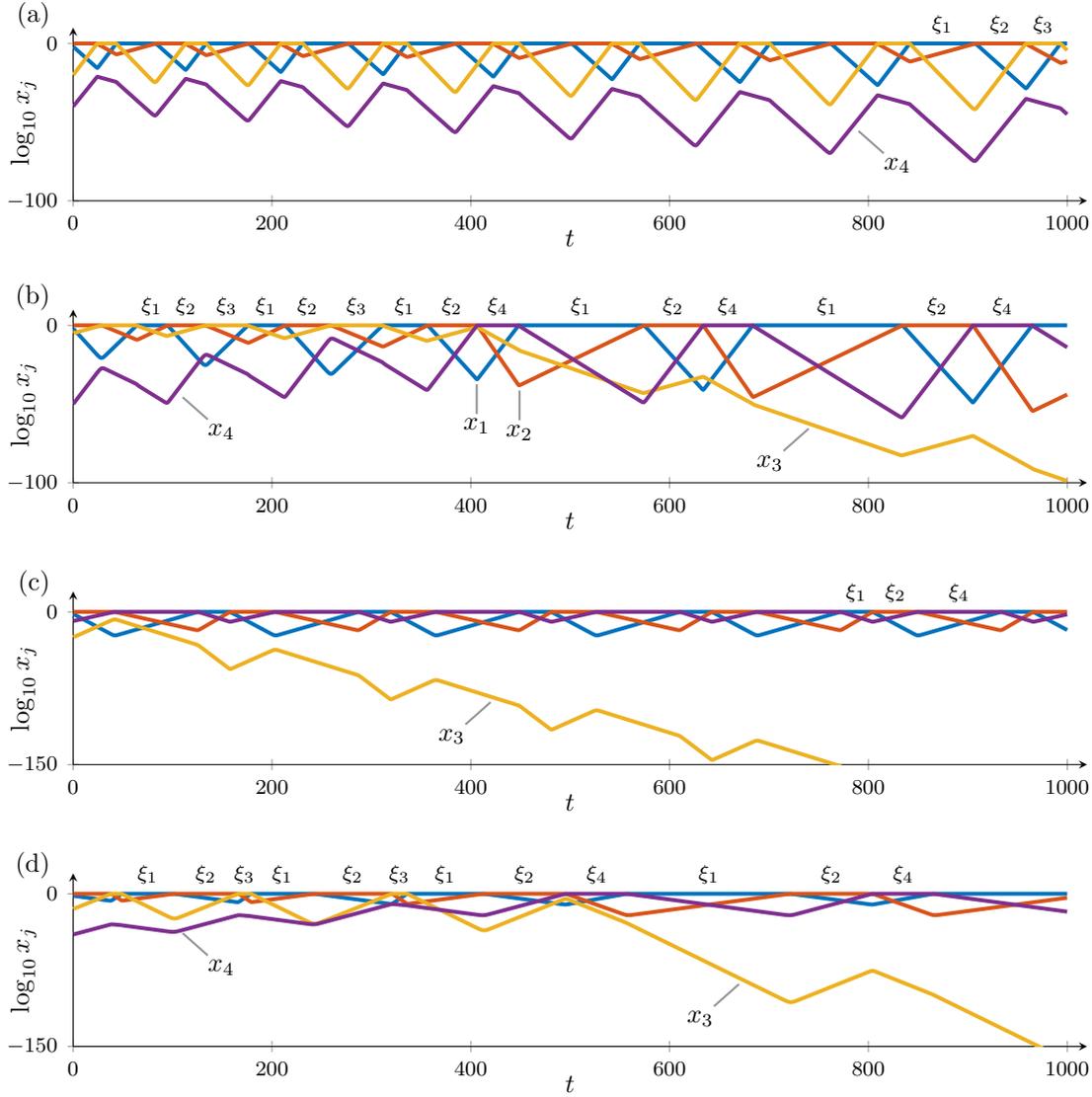
                              
\hbox to \hsize{\hfil%
\fbox{\beginpgfgraphicnamed{x1x2x3x4_ks_timeseries_1}%
\endpgfgraphicnamed}\hfil}          
                          
\hbox to \hsize{\hfil%
\fbox{\beginpgfgraphicnamed{x1x2x3x4_ks_timeseries_2}%
\endpgfgraphicnamed}\hfil}

\hbox to \hsize{\hfil%
\fbox{\beginpgfgraphicnamed{x1x2x3x4_ks_timeseries_3}%
\endpgfgraphicnamed}\hfil}    

\hbox to \hsize{\hfil%
\fbox{\beginpgfgraphicnamed{x1x2x3x4_ks_timeseries_4}%
\endpgfgraphicnamed}\hfil}

\caption{Trajectories near the Kirk--Silber network.
In each case, the coordinates are plotted on a logarithmic axis, each colour corresponding to one coordinate as indicated.
Parameters are given in table~\ref{tab:parametervalues} in appendix A. In (a) (section~\ref{sec:KS1}), the trajectory is approaching the $\xi_1$-$\xi_2$-$\xi_3$ cycle: the time spent near, e.g., $\xi_1$ increases on each loop around the cycle, the minimum value of the $x_1$, $x_2$ and $x_3$ coordinates decreases on each loop, and additionally the $x_4$ coordinate is decaying.
In (b) (section~\ref{sec:KS2}), the trajectory initially appears to approach $\xi_1$-$\xi_2$-$\xi_3$: again the time spent near $\xi_1$ increases on each loop around the cycle.
However, now the $x_4$ coordinate grows and the trajectory switches to the $\xi_1$-$\xi_2$-$\xi_4$ cycle. 
In (c) and (d) the $\xi_1$-$\xi_2$-$\xi_4$ cycle is at resonance.
In (c) (section~\ref{sec:KS3}), the trajectory starts near the $\xi_1$-$\xi_2$-$\xi_4$ cycle and the $x_3$ coordinate decays.
However, trajectories do not get closer to the $\xi_1$-$\xi_2$-$\xi_4$ within that subspace as they remain near a periodic orbit. 
In (d) (section~\ref{sec:KS4}), the trajectory starts near the $\xi_1$-$\xi_2$-$\xi_3$ cycle and then switches to the $\xi_1$-$\xi_2$-$\xi_4$ cycle; the $\xi_1$-$\xi_2$-$\xi_4$ cycle is at resonance so the trajectory remains a finite distance away.
\label{fig:KS_example}}
\end{figure}

We give four examples of parameter choices that illustrate how varied the outcomes can be, thus illustrating the deficiency in the information provided by stating that network is f.a.s. 
Each of these four examples is illustrated by a numerical trajectory shown in one panel of Figure~\ref{fig:KS_example}.
These examples are part of our motivation to introduce finer notions of \emph{visibility} in Section~\ref{sec:visibility}.

\subsubsection{Not at resonance, no switching}
\label{sec:KS1}

For this example, we choose parameters so that $\rho_{123}, \rho_{124}>1$, and all $\nu_{ijk3},\nu_{ijk4}>0$.
Then each of the $\xi_1$-$\xi_2$-$\xi_3$ and $\xi_1$-$\xi_2$-$\xi_4$ subcycles are asymptotically stable for the flow restricted to their own invariant subspaces. Since there is no switching between cycles, then each cycle is f.a.s.\ in the full system. 
Even though the Kirk-Silber network is f.a.s., only one of the two subcycles will appear as the outcome of numerical simulations; which cycle appears depends on the initial conditions.
In Figure~\ref{fig:KS_example}(a) we show a trajectory asymptoting onto the $\xi_1$-$\xi_2$-$\xi_3$ cycle. 

\subsubsection{Not at resonance, with switching}
\label{sec:KS2}

For this example, we choose parameters so that $\rho_{123}, \rho_{124}>1$, all $\nu_{ijk3}>0$, but at least one $\nu_{ijk4}<0$. Then 
each of the $\xi_1$-$\xi_2$-$\xi_3$ and $\xi_1$-$\xi_2$-$\xi_4$ sub-cycles are again asymptotically stable for the flow restricted to their own invariant subspaces. However, now  there is switching from the $\xi_3$ cycle to the $\xi_4$ cycle. 

The $\xi_1$-$\xi_2$-$\xi_4$ cycle is f.a.s., and the $\xi_1$-$\xi_2$-$\xi_3$ cycle is not.  
However, almost all initial conditions near the $\xi_1$-$\xi_2$-$\xi_4$ cycle eventually asymptote onto that cycle, although some may visit the $\xi_1$-$\xi_2$-$\xi_3$ cycle first.  
The $\xi_1$-$\xi_2$-$\xi_4$ cycle will appear in most numerical simulations; even though it is `only' f.a.s., if nearby trajectories move away and make some loops around the $\xi_1$-$\xi_2$-$\xi_3$ cycle, after a transient these trajectories will eventually asymptote onto the $\xi_1$-$\xi_2$-$\xi_4$ cycle.
 
In Figure~\ref{fig:KS_example}(b) we show a trajectory which starts near~$\xi_2$, so the initial conditions is close to both the 
$\xi_1$-$\xi_2$-$\xi_3$ cycle, and the $\xi_1$-$\xi_2$-$\xi_4$ cycle. The trajectory
makes three excursions around the $\xi_1$-$\xi_2$-$\xi_3$ cycle, and then switches onto the $\xi_1$-$\xi_2$-$\xi_4$ cycle.

\subsubsection{At resonance,  no switching}
\label{sec:KS3}

For this example, we choose parameters so that $\rho_{124}=1$, $\rho_{123}>1$, and all $\nu_{ijk3}, \nu_{ijk4}>0$.
Then the  $\xi_1$-$\xi_2$-$\xi_4$ cycle is at resonance, so as with the Guckenheimer--Holmes resonance example (figure~\ref{fig:GHcycnet}(c)), numerical simulations near the $\xi_1$-$\xi_2$-$\xi_4$ cycle shows trajectories approaching a periodic orbit close to that cycle in the subspace $x_3=0$.
The $\xi_1$-$\xi_2$-$\xi_4$ cycle is not f.a.s.  
In Figure~\ref{fig:KS_example}(c) we show a trajectory of this type, and in Figure~\ref{fig:KS_example_linear}(c) we show the same trajectory in linear coordinates. 
The $\xi_1$-$\xi_2$-$\xi_4$ cycle is easily identifiable in the numerics, but unlike in panels (a) and (b), the period of the oscillations does not lengthen.
The dynamics near the $\xi_1$-$\xi_2$-$\xi_3$ cycle are as in Section~\ref{sec:KS1}, that is, the $\xi_1$-$\xi_2$-$\xi_3$ cycle is f.a.s. 
Trajectories which start near the $\xi_1$-$\xi_2$-$\xi_3$ cycle approach that cycle.


\subsubsection{At resonance,  switching}
\label{sec:KS4} 

For our last example, we choose parameters so that $\rho_{124}=1$, $\rho_{123}>1$, all $\nu_{ijk3}>0$, and at least one $\nu_{ijk4}<0$.
So again, the $\xi_1$-$\xi_2$-$\xi_4$ cycle is at resonance and now there is switching from the $\xi_3$ cycle to the $\xi_4$ cycle. 

\begin{figure}
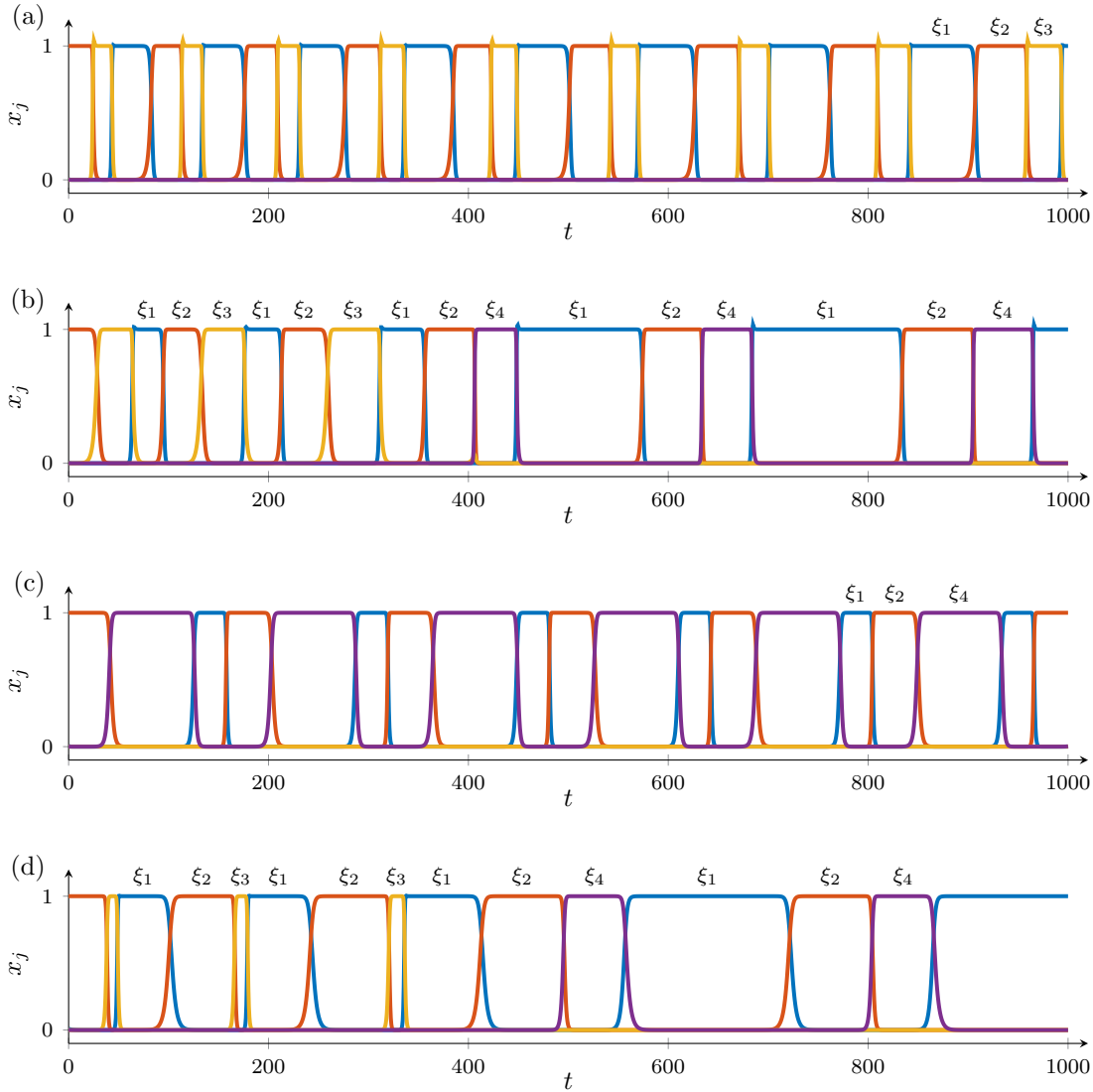
                              
\hbox to \hsize{\hfil%
\fbox{\beginpgfgraphicnamed{x1x2x3x4_ks_timeseries_1_linear}%
\endpgfgraphicnamed}\hfil}          
                          
\hbox to \hsize{\hfil%
\fbox{\beginpgfgraphicnamed{x1x2x3x4_ks_timeseries_2_linear}%
\endpgfgraphicnamed}\hfil}

\hbox to \hsize{\hfil%
\fbox{\beginpgfgraphicnamed{x1x2x3x4_ks_timeseries_3_linear}%
\endpgfgraphicnamed}\hfil}    

\hbox to \hsize{\hfil%
\fbox{\beginpgfgraphicnamed{x1x2x3x4_ks_timeseries_4_linear}%
\endpgfgraphicnamed}\hfil}

\caption{Trajectories near the Kirk--Silber network: the same data as in Figure~\ref{fig:KS_example} but on a linear scale. 
\label{fig:KS_example_linear}}
\end{figure}

Trajectories which start near the $\xi_1$-$\xi_2$-$\xi_4$ cycle eventually approach a periodic orbit close to that cycle in the subspace $x_3=0$.
Some trajectories (specifically, an open set starting near $\xi_2$) will first make one, or several, excursions around the $\xi_1$-$\xi_2$-$\xi_3$ cycle before switching back to the $\xi_1$-$\xi_2$-$\xi_4$ cycle. 
Thus the $\xi_1$-$\xi_2$-$\xi_4$ cycle does not satisfy the definitions of either Lyapunov stability or quasi-asymptotic stability.
However, trajectories near the $\xi_1$-$\xi_2$-$\xi_4$ cycle will be observable in numerical simulations.
Figure~\ref{fig:KS_example}(d) shows such a trajectory. The trajectory is shown again in linear coordinates in Figure~\ref{fig:KS_example_linear}(d). It can be seen here that the period of oscillations as the trajectory goes around $\xi_1$-$\xi_2$-$\xi_3$ is increasing, but after the switch to $\xi_1$-$\xi_2$-$\xi_4$, the period remains constant.
The $\xi_1$-$\xi_2$-$\xi_3$ cycle also does not satisfy any definition of stability as all trajectories which start near it will eventually switch to be near the $\xi_1$-$\xi_2$-$\xi_4$ cycle.

\subsection{Rock--Paper--Scissors--Spock--Lizard}\label{sec:LS}

The Rock--Paper--Scissors--Spock--Lizard network is a union of five Guckenheimer--Holmes cycles and is shown in Figure~\ref{fig:KS_and_RPSSL_networks}(b).
The dynamics near this network are highly complex and some results on the stability of subcycles and omnicycles can be found in~\cite{Postlethwaite2022} and~\cite{Castro2022a} (see~\cite{Podvigina2023} for the definition of omnicycle).
Here we give four example numerical simulations where different parts of the network, but not necessarily the whole network, appear in the long term. 
In each case, we are in a parameter region where the entire network has been shown to be asymptotically stable.

The network exists in the following equations:
\begin{align}
\dot{x}_1& = x_1\left(1-\chi-c_A x_2^2+ e_B x_3^2 -c_B x_4^2+ e_A x_5^2\right) \nonumber \\
\dot{x}_2& = x_2\left(1-\chi-c_A x_3^2+ e_B x_4^2 -c_B x_5^2+ e_A x_1^2\right) \nonumber \\
\dot{x}_3 & = x_3\left(1-\chi-c_A x_4^2+ e_B x_5^2 -c_B x_1^2+ e_A x_2^2\right) \label{eq:odes} \\
\dot{x}_4 & = x_4\left(1-\chi-c_A x_5^2+ e_B x_1^2 -c_B x_2^2+ e_A x_3^2\right) \nonumber \\
\dot{x}_5 & = x_5\left(1-\chi-c_A x_1^2+ e_B x_2^2 -c_B x_3^2+ e_A x_4^2\right) \nonumber
\end{align}
where again $\chi=\sum_{j=1}^5{x_j^2}$ and $e_A,e_B,c_A,c_B>0$ are parameters. As with our other examples, there are equilibria lying on each coordinate axis, each with one non-zero coordinates. We refer to these equilibria as $\xi_j$ ($j=1,\dots, 5$).

In~\cite{Postlethwaite2022} the edges of the Rock--Paper--Scissors--Spock--Lizard network are divided into two types: A, which are heteroclinic connections between equilibria $\xi_j$ and $\xi_{j+1}$, and B, which are heteroclinic connections between equilibria $\xi_j$ and $\xi_{j-2}$ (with indices taken mod $5$). 
The B-type connections  are one-dimensional, but the $A$ connections are actually a two-dimensional manifold of connecting orbits. 
The boundary of this two-dimensional manifold consists of a one-dimensional connecting trajectory from $\xi_j$ to $\xi_{j+1}$, and two B-connections from $\xi_j$ to $\xi_{j-2}$ and $\xi_{j-2}$ to $\xi_{j-4}\equiv \xi_{j+1}$.
{Given our definition of heteroclinic cycle, when confronted with higher-dimensional connections, we select the one-dimensional connecting trajectory on the boundary to construct the heteroclinic cycle.
This trajectory is contained in a coordinate plane.
We still use type~A when referring to it and explicitly provide the dimension when necessary.}

\begin{figure}
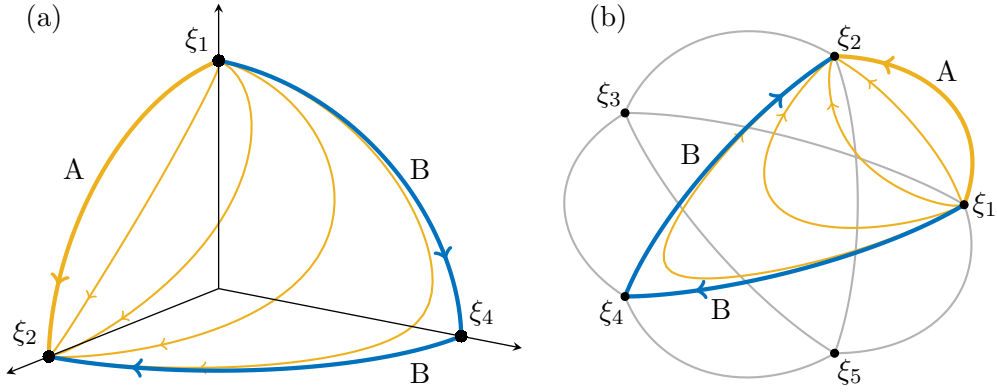

\hbox to \hsize{\hfil%
\fbox{\beginpgfgraphicnamed{RPSSL_1D_2D_manifolds}%
\endpgfgraphicnamed}\hfil}
\caption{The figures show the two-dimensional unstable manifold of $\xi_1$ in the RPSSL network. In (a), a set of trajectories within this manifold is shown in the three-dimensional $x_1$-$x_2$-$x_4$-subspace. In (b), the same trajectories are shown using the pentacle projection described in the text (the origin is in the center of the figure, and $\xi_1$ is at $(1,0)$). 
A~and B type connections are labelled in each figure. 
The amber trajectories are part of a two-dimensional connection of type~A.
The connecting trajectory in the plane containing the equilibria $\xi_1$ and $\xi_2$ is the one-dimensional representative trajectory of type~A. 
The choice of colours is the same as in Figure~\ref{fig:KS_and_RPSSL_networks}(b).
The light gray lines show the remaining one-dimensional A and B connections (compare with figure~\ref{fig:KS_and_RPSSL_networks}(b)).
 \label{fig:2Dwu}}
\end{figure}


For the clarity of displaying our numerical results, we use what we term the \emph{pentacle projection}~\cite{Ashwin2013} to project the five-dimensional phase space onto two dimensions. Consider a trajectory $\mathbf{x}(t)=(x_1(t),x_2(t),x_3(t),x_4(t),x_5(t))\in\R^5$. The pentacle projection is the map:
\[
P:\R^5\rightarrow \R^2 
\]
where we write
\[
y_1(t)=\sum_{j=1}^5 x_j(t) \cos\left( \frac{2\pi j}{5}\right), \qquad  y_2(t)=\sum_{j=1}^5 x_j(t) \sin\left( \frac{2\pi j}{5}\right),
\]
and then $P(\mathbf{x})=\mathbf{y}\equiv(y_1(t),y_2(t))\in\R^2$.
Notice that the points $P(\xi_1),\dots,P(\xi_5)$ lie on the vertices of a pentagon. In Figure~\ref{fig:2Dwu} we show the two-dimensional unstable manifold of $\xi_1$ both in a three-dimensional subspace and also in the pentacle projection.

\begin{figure}[htp]
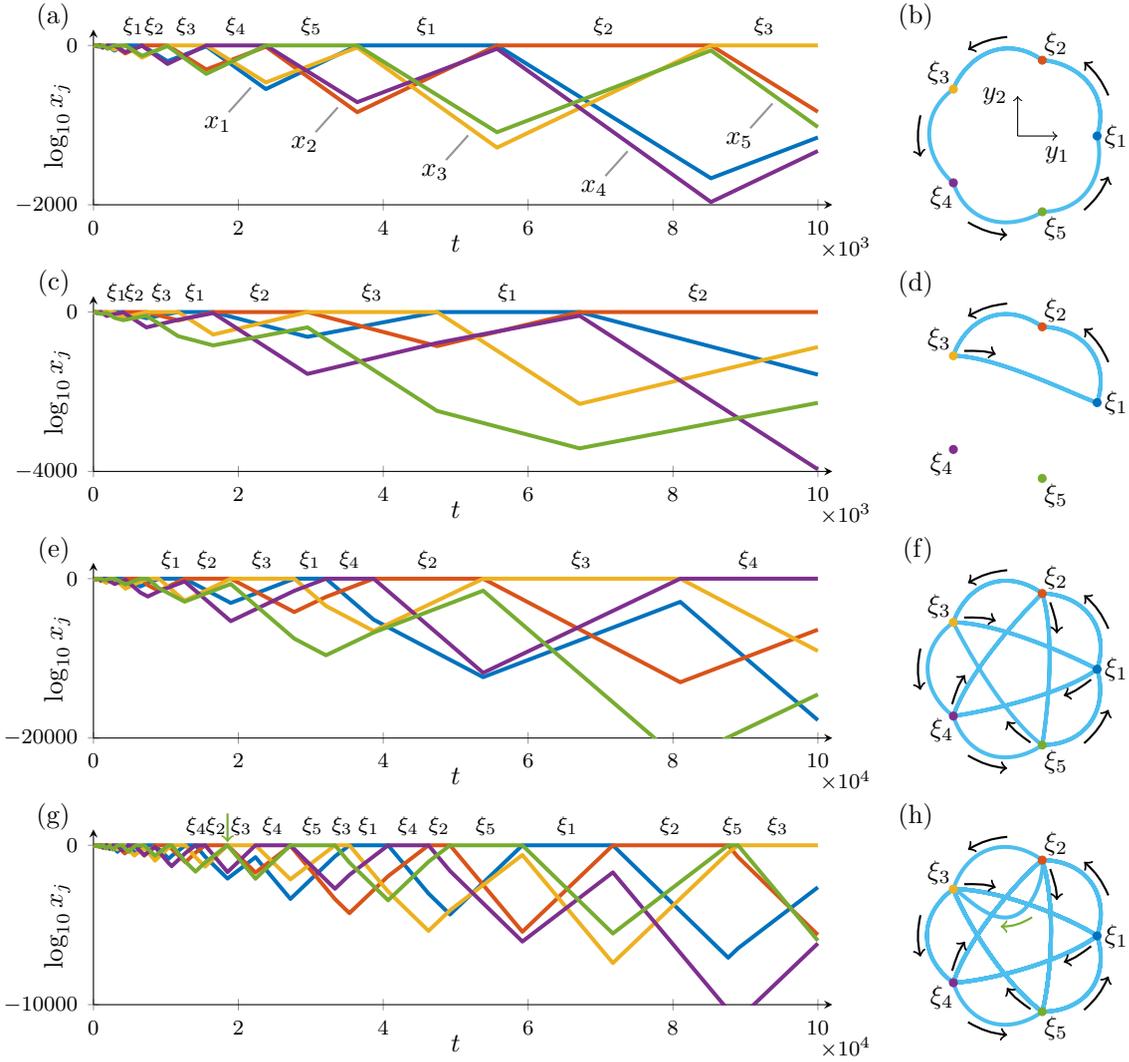

\hbox to \hsize{\hfil%
\fbox{\beginpgfgraphicnamed{x1x2x3x4x5_rpssl_timeseries_1}%
\endpgfgraphicnamed}\hfil
\fbox{\beginpgfgraphicnamed{x1x2x3x4x5_rpssl_pentacle_1}%
\endpgfgraphicnamed}\hfil}

\hbox to \hsize{\hfil%
\fbox{\beginpgfgraphicnamed{x1x2x3x4x5_rpssl_timeseries_2}%
\endpgfgraphicnamed}\hfil
\fbox{\beginpgfgraphicnamed{x1x2x3x4x5_rpssl_pentacle_2}%
\endpgfgraphicnamed}\hfil}

\hbox to \hsize{\hfil%
\fbox{\beginpgfgraphicnamed{x1x2x3x4x5_rpssl_timeseries_3}%
\endpgfgraphicnamed}\hfil
\fbox{\beginpgfgraphicnamed{x1x2x3x4x5_rpssl_pentacle_3}%
\endpgfgraphicnamed}\hfil}

\hbox to \hsize{\hfil%
\fbox{\beginpgfgraphicnamed{x1x2x3x4x5_rpssl_timeseries_4}%
\endpgfgraphicnamed}\hfil
\fbox{\beginpgfgraphicnamed{x1x2x3x4x5_rpssl_pentacle_4}%
\endpgfgraphicnamed}\hfil}

\caption{In each pair of panels, one trajectory is shown: on the left as a time series in logarithmic coordinates, and on the right in the pentacle projection. In (a) and (b) (section~\ref{sec:RPSSLA}) the trajectory approaches the $A$-subcycle. In (c) and (d) (section~\ref{sec:RPSSLAAB}) the trajectory approaches one of the $AAB$-cycles, in this case, the $\xi_1$-$\xi_2$-$\xi_3$ cycle. In (e) and (f) (section~\ref{sec:RPSSLAABBB}), the $AABBB$ omnicycle is f.a.s., and the trajectory visits the entire {one-dimensional} RPSSL network. In (g) and (h) (section~\ref{sec:RPSSLchaos}), the trajectory shown visits the equilibria in an aperiodic fashion. The entire {one-dimensional} RPSSL network is visited, but in addition the trajectory sometimes makes excursions through the interior of the two-dimensional unstable manifolds, indicated by the green arrow between $\xi_2$ and~$\xi_3$ in (g) and~(h). In all cases, the panels in the right column have the transient up to $t=500$ removed. Parameters are given in table~\ref{tab:parametervalues} in appendix A.
\label{fig:RPSSL_example}}
\end{figure}

\subsubsection{$A$ cycle}
\label{sec:RPSSLA}
Here we choose parameters so that the $A$-subcycle, {consisting of one-dimensional connecting trajectories on the coordinate planes,} is f.a.s.
In Figure~\ref{fig:RPSSL_example}(a) and (b) we show a trajectory which asymptotes onto this subcycle ($\xi_1$-$\xi_2$-$\xi_3$-$\xi_4$-$\xi_5$). Although it has not been proven that there are no other attracting dynamics for these parameter value, extensive numerical investigations have been unable to find any. Numerical experiments strongly suggest that almost all initial conditions asymptote onto the $A$-cycle, and as such, only the $A$-type edges of the network are seen in the long term.

 \subsubsection{$AAB$ cycle}
\label{sec:RPSSLAAB}

Here we choose parameters so that each $AAB$ subcycle (of which there are five) is f.a.s. 
Note that the $AAB$ subcycle, {when considering one-dimensional connecting trajectories in coordinate planes,} is equivalent to a Guckenheimer and Holmes cycle (with a slight restriction on the choice of parameter values).
The cycle that appears in numerical simulations will depend on the initial conditions. In Figure~\ref{fig:RPSSL_example}(c) and (d) we show a trajectory that asymptotes onto the $\xi_1$-$\xi_2$-$\xi_3$ subcycle. 
Again, it has not been proven that there are no other attracting dynamics for these parameter values but numerical experiments strongly suggest this is the case. 
Trajectories will not explore the entire network in the long term, but will be attracted to a sub-cycle of length three.

 \subsubsection{$AABBB$ cycling}
 \label{sec:RPSSLAABBB}

Here we choose parameter such that the $AABBB$ omnicycle is f.a.s. That is, there is an open set of initial conditions such that trajectories repeated follow the edges in the order $AABBB$ and also get closer to the network. 
A trajectory of this type is shown in Figure~\ref{fig:RPSSL_example}(e) and (f). 
In this case, the entire {one-dimensional} RPSSL network is visited by the trajectory. 

\subsubsection{Aperiodic behaviour}
\label{sec:RPSSLchaos}

Our final example in this section uses parameters for which the $A$-cycle is f.a.s.
In addition, numerical experiments have found trajectories that perform apparently irregular cycling around the network: the observed sequence of edges $A$ and $B$ that is followed by the trajectory is aperiodic. 
An example of such as trajectory is shown in Figure~\ref{fig:RPSSL_example}(g) and (h). 
The entire network of one-dimensional connections is visited, but in addition, the trajectory occasionally moves away from this set of one-dimensional connections and passes through the interior of the two-dimensional connection from $\xi_j$ to $\xi_{j+1}$. 
This can be seen in Figure~\ref{fig:RPSSL_example}(h) where there is a curve from $\xi_2$ to $\xi_3$ that is not close to the {one-dimensional} RPSSL network. 

We suspect, but have not proved, that for a long enough trajectory, the entirety of the two-dimensional RPSSL network would be visited by the trajectory, in the sense that the trajectory would eventually come arbitrarily close to any point on the two-dimensional manifolds.

\section{Visibility}\label{sec:visibility}

The examples in Section~\ref{sec:examples} clearly indicate that establishing the stability properties of a heteroclinic network and its cycles does not always provide information about what to expect from numerical simulations, especially in the long-run.
In this section, without attempting to cover all possible instances, we suggest the use of the notion of \emph{visibility} to express the likely outcome of the dynamics near a heteroclinic network.
We are particularly concerned with the fact that a network can be called `asymptotically stable' (or a weaker version of this) while trajectories could visit only a part of the network, with the remaining part never seen in numerical simulations, or only seen during a transient.
We introduce the concept of `visibility' to mean a set that is attracting and nearby trajectories repeatedly visit the \emph{entire} set.

For example, in the context of population dynamics when the concern is to establish whether various species survive, the stability of the entire network does not give this information. 
Take the case of five species competing as in the RPSSL. In all four of the examples we give above, the entire RPSSL network is asymptotically stable, but this does not establish the survival of all the species. 
For example in the second example (Section~\ref{sec:RPSSLAAB}), only three of the species represented survive, and which ones depends on the initial conditions. 
The concept of visibility will classify each of the five three-species subcycles as visible, but not the whole network.




\subsection{Definitions of visibility}

We align our definitions of visibility (and all its variants) along the lines of the previous definitions of stability. In particular:
\begin{enumerate}[label=(\roman*)]

\item  we maintain the use of prefixes `fragmentary-' or `essential-' when describing visible invariant sets to clarify \emph{from where} the set is visible, that is, from which initial conditions near the set.
If neither are used the invariant set is visible from a neighbourhood of the set. 
In addition, we introduce the prefix `almost-' to mean a set is visible from a neighbourhood, excluding a set of measure zero.
    
\item  we maintain the use of the prefix `quasi-' when trajectories may contain an initial transient which moves away from the visible invariant set. 
That is, we allow visible sets to not be visible for an initial transient.
This clarifies \emph{from when} the set is visible.

\item we maintain the use of the prefix `Lyapunov' to indicate that trajectories should `stay close' to the visible set.
\end{enumerate}

In parallel with the corresponding definitions for stability, we define \emph{Lyapunov-visibility} to mean sets for which a neighbourhood of initial conditions stay close to the set -- and visit the entire set (but may not asymptote onto that set), and \emph{quasi-visibility} to mean sets for which a neighbourhood of initial conditions eventually asymptote onto the entire set (but may move a large distance away from the set initially). 
Similarly, if a set is both Lyapunov- and quasi-visible, then we say it is \emph{asymptotically visible}. 
The precise definitions we propose are given next.

Let $N_{\delta}=\{x \notin X$ : $d(x,X)<\delta\}$. This set is a ball around $X$ excluding $X$ itself. We do this because  if $X$ is a heteroclinic cycle or network, points in $X$ (that is, equilibria and connecting trajectories), never visit the whole of $X$.

\begin{definition}\label{def:L-visible}
    A closed invariant set $X$ is called \emph{Lyapunov-visible}, or \emph{L-visible}, if for any~$\varepsilon>0$ there is a~$\delta>0$ such that for all $\bar{x} \in X$ and for all $x \in N_{\delta}$, $d(\varphi(t,x),X)<\varepsilon$ for all $t>0$, and there exists an increasing sequence $(t_n)$ with $\lim_{n\rightarrow\infty} t_n =\infty$ such that $d(\varphi(t_n,x),\bar{x})<\varepsilon$.
\end{definition}
Demanding that trajectories repeatedly come within $\varepsilon$ of every point in $X$ ensures that trajectories starting at initial conditions near $X$ visit the whole of $X$ (and not just a subset of $X$).  

\begin{definition}\label{def:Q-visible}
A closed invariant set $X$ is called \emph{quasi-visible}, or \emph{Q-visible},  if there exists $\delta>0$ such that if an initial condition
$x \in N_{\delta}$ then $\omega(x)=X$. 
\end{definition}
Notice that this definition differs from Definition~\ref{def:quasi-stability} by requiring that $\omega(x)=X$, so trajectories must visit the whole of $X$.

\begin{definition}\label{def:As-visible}
    A closed invariant set $X$ is called \emph{asymptotically visible}, or \emph{A-visible} if it is both Q-visible and L-visible.
\end{definition}

\begin{definition}\label{def:visible}
    A closed invariant set $X$ is called \emph{visible} if for any~$\varepsilon>0$ there is a~$\delta>0$ such that for all $\bar{x} \in X$, and for all $x \in N_{\delta}$,  there exists a $T>0$, such that $d(\varphi(t,x),X)<\varepsilon$ for all $t>T$, and there exists an increasing sequence $(t_n)$ with $\lim_{n\rightarrow\infty} t_n =\infty$ such that $d(\varphi(t_n,x),\bar{x})<\varepsilon$.
\end{definition}

The difference between Definitions~\ref{def:L-visible} and~\ref{def:visible} is that in~\ref{def:visible}, we remove a transient, whose length may depend on the initial condition, before insisting that the trajectory remains close to the entire set~$X$.
It is important to note that a visible set could be neither Q- nor L-visible. 
With this definition it is possible to have some nearby points move away first, and come back, but no points ever asymptote onto the set.

In Definitions~\ref{def:L-visible}--\ref{def:visible} we use neighbourhoods of~$X$ excluding points within~$X$ itself. 
Thus a set $X$ is Q-visible if it eventually attracts a neighbourhood of points, not in~$X$, which are allowed to diverge from~$X$ for a finite time. 
Analogously, L-visible sets are such that initial conditions in a neighbourhood of~$X$, but not in~$X$, do not move far away.
When studying heteroclinic cycles in a network the use of neighbourhoods is excessively strong even when~$X$ itself is excluded, for two reasons.
Firstly, systems containing robust heteroclinic connections typically contain invariant subspaces. 
Trajectories that start in these subspaces will necessarily not visit the entire network (due to the invariance). 
These subspaces will have measure zero in the phase space. 
Secondly, a neighbourhood of an equilibrium with more than one outgoing direction necessarily has points following both heteroclinic connections, and thus requiring all trajectories in a neighbourhood to follow a particular path is overly restrictive.
In order to include heteroclinic cycles in networks in our study, we use the prefixes almost-, essentially- and fragmentarily-  to allow for the exclusion of sets (of various sizes) from a neighbourhood of $X$. 


\begin{definition}\label{def:almost-visible}
     A closed invariant set $X$ is called \emph{almost} Lyapunov-/quasi-/a\-symp\-to\-tically-visible, if the set $N_{\delta}$ in the definitions for Lyapunov-/quasi-/asymptotic- visibility is replaced by a set $M_{\delta}$ which differs from $N_{\delta}$ by only a set of measure zero.
\end{definition}


\begin{definition}\label{def:E-visible}
     A closed invariant set $X$ is called \emph{essentially} Lyapunov-/quasi-/a\-symp\-to\-tically- visible, if the set $N_{\delta}$ in the definitions for Lyapunov-/quasi-/asymptotic- visibility is replaced by a set $M_{\delta}\subset N_{\delta}$ that satisfies
    \begin{align*}
        \lim\limits_{\varepsilon \to 0} \frac{\ell(V_{\varepsilon}(X) \cap M_{\delta})}{\ell(V_\varepsilon(X))} = 1,
    \end{align*}
    where $V_{\varepsilon}$ is an $\varepsilon$-neighbourhood of~$X$.
\end{definition}


\begin{definition}\label{def:F-visible}
     A closed invariant set $X$ is called \emph{fragmentarily} Lyapunov-/qua\-si-/asymp\-totically visible, if 
     the set $N_{\delta}$ in the definitions for Lyapunov-/quasi-/asymptotic- visibility is replaced by a set $M_{\delta}\subset N_{\delta}$ that has positive measure.
\end{definition}

Because of the size of $M_\delta$ in each of these definitions, `almost-' implies \mbox{`essentially-'}\ implies `fragmentarily-'.
Note that we are using the prefixes fragmentarily- and essentially- \emph{without} requiring the Lyapunov prefix. 
This is different to the original definition of f.a.s.~\cite{Podvigina2012}. 
For instance, in Example~\ref{sec:KS2}, the $\xi_1$-$\xi_2$-$\xi_4$ cycle is fragmentarily-quasi-visible. 
A positive measure set of points which start in a neighbourhood of that cycle eventually asymptote onto that cycle, but some of these trajectories will first make an excursion around the $\xi_1$-$\xi_2$-$\xi_3$ cycle. 
Kirk and Silber~\cite{Kirk1994} refer to this behaviour as `asymptotically stable in spirit'. 

\subsection{Application of visibility definitions to examples}

The Guckenheimer--Holmes cycle in Section~\ref{sec:GH} is asymptotically-visible when $\rho_{123}>1$ and Lyapunov visible when $\rho_{123}=1$.
In this example there is a correspondence between stability and visibility due to the simplicity of the system: there is only one route trajectories can take around the cycle.

For the Kirk--Silber network, we note that it is never possible for the entire network to be visible: trajectories only ever asymptote onto one sub-cycle of the network. However, the sub-cycles can individually be visible.
For the example in Subsection~\ref{sec:KS1}, both cycles are fragmentarily asymptotically-visible. 
In Subsection~\ref{sec:KS2}, 
the $\xi_1$-$\xi_2$-$\xi_4$ cycle is quasi-visible; the $\xi_1$-$\xi_2$-$\xi_3$ cycle is not visible. 
In Subsection~\ref{sec:KS3}, the $\xi_1$-$\xi_2$-$\xi_3$ cycle is fragmentarily asymptotically-visible and the $\xi_1$-$\xi_2$-$\xi_4$ cycle is fragmentarily Lyapunov-visible. 
Finally, in Subsection~\ref{sec:KS4}, the $\xi_1$-$\xi_2$-$\xi_3$ cycle is not visible, but the $\xi_1$-$\xi_2$-$\xi_4$ cycle is fragmentarily visible. 
Here the $\xi_1$-$\xi_2$-$\xi_4$ cycle is neither fragmentarily-quasi-visible, because trajectories remain a finite distance away from the cycle due to the resonance condition; nor is it fragmentarily-Lyapunov-visible, because some nearby trajectories will make one or more excursions around the $\xi_1$-$\xi_2$-$\xi_3$ cycle first.

The {one-dimensional} RPSSL  network, in contrast to the KS network, can be visible, depending on parameters. 
In Subsection~\ref{sec:RPSSLA}, results from~\cite{Postlethwaite2022,Castro2022a} prove that the $A$-cycle is fragmentarily-asymptotically visible. 
Our numerical experiments strongly suggest that the entire network is not visible, and in fact that the $A$-cycle has stronger visibility: that it is asymptotically visible. 
In Subsection~\ref{sec:RPSSLAAB}, each of the five $AAB$-cycles is fragmentarily-asymptotically visible. 
Again, numerical experiments strongly suggest that the entire network is not visible. 
In Subsection~\ref{sec:RPSSLAABBB}, results from~\cite{Postlethwaite2022} can be used to prove the entire network is fragmentarily visible, but again we conjecture that in fact the network is asymptotically visible. 
In Subsection~\ref{sec:RPSSLchaos} the numerical experiments suggest the network is not visible. Although the irregular trajectories repeated visit all parts of the network, there are rare excursions away from the {one-dimensional} network. 
We conjecture that these excursions continue to occur indefinitely, thus the trajectories do not satisfy any definitions of visibility. 
Furthermore, we conjecture that the 2D RPSSL network is visible, as eventually all parts of the two-dimensional unstable manifold will be visited, but a proof of this is beyond the scope of this paper.


\section{Discussion}

In this paper we have introduced a new concept, visibility, which we believe is as important as the stability of an invariant object when one is concerned with predicting observations of a dynamical system after transients have decayed. 
Our examples are all heteroclinic cycles or networks; it has been clear for some years that traditional stability measures such as asymptotic stability are not necessarily the most useful in these cases, which is why weaker properties such as fragmentary asymptotic stability have been introduced. 
To determine the visibility of an object one must asks additional questions to determine whether the entire object will be visited by trajectories after transients have decayed.

This is of course not to say that stability is unimportant or not useful. 
For example, in game theory the outcome of a game is frequently a Nash equilibrium. In a class of games known as `coordination games' there is more than one Nash equilibrium and deciding which one occurs is still an open and relevant question. 
Recently, Castro~\cite{Castro2025} has used techniques from heteroclinic dynamics to, under some assumptions, decide which Nash equilibrium is chosen (in this case, one that is asymptotically stable in an extension of the original game). 

In this paper we have not provided any proofs, but we note that proofs of some of our statements can be deduced from previously known results. 
For example, results in the original paper by Kirk and Silber~\cite{Kirk1994} show that the Kirk--Silber network cannot be visible, for any choices of parameters.
Results from~\cite{Postlethwaite2022} showing that the $AABBB$ cycle of the RPSSL network is f.a.s.\ imply that the network is fragmentarily-asymptotically visible.

We conclude by offering several avenues of work which require further investigation. 
First, we note that proving whether or not a network is visible will typically require some graph theory as well as stability calculations. 
For example, in the RPSSL network, stabililty of any of the omnicycles $A$, $B$ or $AAB$ does \emph{not} imply visibility of the whole network, but stability of $AABBB$ does. 
This is due to the structure of the digraph representing the network and whether repeating a certain sequence of types of edges generates a path which covers the entire graph. 
For larger networks/digraphs, this could become computationally difficult.

Second, the issue of two-dimensional unstable manifolds requires special care. 
Our definition of heteroclinic networks in the introduction implies that they are one-dimensional objects. 
However, many heteroclinic networks have at least one equilibria which has an unstable manifold of dimension two or higher. 
In some cases, such as the KS-network, it makes sense to study the visibility, or stability of the one-dimensional set of heteroclinic connections which intersect with the invariant coordinate planes. 
However, our numerical computations near the RPSSL network, as well as in some other larger networks, have indicated that there exists aperiodic switching where trajectories repeatedly move away from this one-dimensional set of connections and pass through the interior of the two-dimensional unstable manifold. 
In this case, the original (one-dimensional) network does not satisfy our definition of visible as trajectories repeatedly move away from it. 
We might ask the question of whether we could include the entire two-dimensional manifold in the invariant object we study, but then a proof of visibility would require proving that trajectories repeatedly visit the entire two-dimensional manifold. How to construct such a proof is beyond our current understanding.

Some of these issues may be able to be better understood using techniques recently introduced to examine Poincar\'{e} maps near heteroclinic networks using techniques from piecewise dynamical systems~\cite{GroothuizenDijkema2025}. 
These methods allows us both to have some level of control over the description of the observed transients and also potentially the opportunity to prove statments about the existence of chaotic attractors.

\section*{Acknowledgments}
We would like to acknowledge conversations with Christian Bick and Alexander Lohse.

SC was partially supported by CMUP, member of LASI, which is financed by national funds through FCT -- Funda\c{c}\~ao para a Ci\^encia e a Tecnologia, I.P., under the projects UIDB/00144/2020 and UIDP/00144/2020. 
The work of CMP was supported by the Marsden Fund Council from New Zealand Government funding, managed by the Royal Society Te Ap\={a}rangi (Grant No.\ 21-UOA-048). 
AMR is supported by the EPSRC grant EP/V014439/1, which is gratefully acknowledged.
AMR and CMP are grateful to the London Mathematical Society for a `Visits to the UK -- Scheme~2' grant.

The data associated with this paper are openly available from the University of Leeds Data Repository (\url{http://doi.org/10.5518/1652})~\cite{Castro2025b}, as are the programs that generated the data.
For the purpose of open access, the authors have applied a Creative Commons Attribution (CC~BY) license to any Author Accepted Manuscript version arising from this submission.

\appendix
\section{Parameters used in figures}

Table~\ref{tab:parametervalues} gives the parameters used in the creation of
figures~\ref{fig:GH_example},~\ref{fig:KS_example},~\ref{fig:KS_example_linear} and~\ref{fig:RPSSL_example}.

\begin{table}
\caption{Parameter values for the Guckenheimer--Holmes example 
(Figure~\ref{fig:GH_example}), the Kirk--Silber examples (Figures~\ref{fig:KS_example} and~\ref{fig:KS_example_linear})
and the Rock--Paper--Scissors--Spock--Lizard examples (Figure~\ref{fig:RPSSL_example}).
 \label{tab:parametervalues}}
 \begin{center}\begin{tabular}{|c|l|l}
 Cycle and figure & Parameter values\\
 \noalign{\vspace{4pt}}
 \hline
 \noalign{\vspace{4pt}}
 \makecell{Guckenheimer--\\
           --Holmes\\[0.5ex]
           Figure~\ref{fig:GH_example}} &
 \makecell{$e_{12}=0.9$, $c_{13}=1.0$.\\
           $e_{23}=1.5$, $c_{21}=0.9$.\\
           $e_{31}=0.6$, $c_{32}=1.2$.}\\
 \noalign{\vspace{4pt}}
 \hline
 \noalign{\vspace{4pt}}
 \makecell{Kirk--Silber\\[0.5ex]
           Figures~\ref{fig:KS_example} and~\ref{fig:KS_example_linear}} &
 \makecell{(a) $e_{12}=0.4$,
               $c_{13}=1.5$,
               $c_{14}=1.3$,
               $c_{21}=1.3$,
               $e_{23}=1.9$,
               $e_{24}=1.8$,\phantom{00}\\
 \phantom{(a)} $e_{31}=1.9$,
               $c_{32}=0.8$,
               $c_{34}=0.4$,
               $e_{41}=1.8$,
               $c_{42}=0.8$,
               $c_{43}=1.2$.\phantom{00}\\
           (b) $e_{12}=0.7$,
               $c_{13}=0.5$,
               $c_{14}=0.9$,
               $c_{21}=1.6$,
               $e_{23}=0.4$,
               $e_{24}=1.9$,\phantom{00}\\
 \phantom{(b)} $e_{31}=1.4$,
               $c_{32}=0.6$,
               $c_{34}=0.7$,
               $e_{41}=1.9$,
               $c_{42}=2.1$,
               $c_{43}=0.8$.\phantom{00}\\
           (c) $e_{12}=1.3$,
               $c_{13}=1.7$,
               $c_{14}=0.7$,
               $c_{21}=1.2$,
               $e_{23}=1.0$,
               $e_{24}=0.5$,\phantom{00}\\
 \phantom{(c)} $e_{31}=2.0$,
               $c_{32}=1.4$,
               $c_{34}=0.5$,
               $e_{41}=0.646$,
               $c_{42}=0.5$,
               $c_{43}=0.7$.\\
           (d) $e_{12}=0.3$,
               $c_{13}=1.1$,
               $c_{14}=0.3$,
               $c_{21}=0.3$,
               $e_{23}=0.9$,
               $e_{24}=0.6$,\phantom{00}\\
 \phantom{(d)} $e_{31}=1.5$,
               $c_{32}=1.5$,
               $c_{34}=0.2$,
               $e_{41}=0.4$,
               $c_{42}=0.8$,
               $c_{43}=0.9$.\phantom{00}}\\
 \noalign{\vspace{4pt}}    
 \hline
 \noalign{\vspace{4pt}}
 \makecell{Rock--Paper--\\
           --Scissors--Spock--\\[0.5ex]
           --Lizard\\[0.5ex]
           Figure~\ref{fig:RPSSL_example}} &
       \makecell{%
       (a) $c_A = 1.30$, $e_A = 1.00$, $c_B = 1.50$, $e_B = 0.80$.\\
       (b) $c_A = 1.10$, $e_A = 1.00$, $c_B = 2.70$, $e_B = 0.80$.\\
       (c) $c_A = 1.10$, $e_A = 1.00$, $c_B = 1.80$, $e_B = 0.80$.\\
       (d) $c_A = 1.02$, $e_A = 1.00$, $c_B = 1.25$, $e_B = 0.80$.}\\
 \noalign{\vspace{4pt}}
 \hline
 \end{tabular}\end{center}
 \end{table}


\begin{thebibliography}{10}

\bibitem{Aguiar2006}
{\sc M.~A.~D. Aguiar, S.~B. S.~D. Castro, and I.~S. Labouriau}, {\em Simple
  vector fields with complex behavior}, Int. J. Bifurcation Chaos, 16 (2006),
  pp.~369--381, \url{https://doi.org/10.1142/S021812740601485X}.

\bibitem{Armbruster1991}
{\sc D.~Armbruster and P.~Chossat}, {\em Heteroclinic orbits in a spherically
  invariant system}, Physica D, 50 (1991), pp.~155--176,
  \url{https://doi.org/10.1016/0167-2789(91)90173-7}.

\bibitem{Ashwin1997f}
{\sc P.~Ashwin}, {\em Cycles homoclinic to chaotic sets; robustness and
  resonance}, Chaos, 7 (1997), pp.~207--220,
  \url{https://doi.org/10.1063/1.166221}.

\bibitem{Ashwin2013}
{\sc P.~Ashwin and C.~Postlethwaite}, {\em {On designing heteroclinic networks
  from graphs}}, {Physica D}, {265} ({2013}), pp.~{26--39},
  \url{https://doi.org/10.1016/j.physd.2013.09.006}.

\bibitem{Auslander1964}
{\sc J.~Auslander, N.~P. Bhatia, and P.~Seibert}, {\em Attractors in dynamical
  systems}, Bol. Soc. Mat. Mex. (2),  (1964), pp.~55--66,
  \url{https://boletin.math.org.mx/pdf/2/9/BSMM(2).9.55-66.pdf}.

\bibitem{Brannath1994}
{\sc W.~Brannath}, {\em Heteroclinic networks on the tetrahedron},
  Nonlinearity, 7 (1994), pp.~1367--1384,
  \url{https://doi.org/10.1088/0951-7715/7/5/006}.

\bibitem{Castro2025}
{\sc S.~B. S.~D. Castro}, {\em Learning coordination through new actions}, Int.
  Game Theory Rev., to appear (2025),
  \url{https://doi.org/10.1142/S0219198924500154}.

\bibitem{Castro2022a}
{\sc S.~B. S.~D. Castro, A.~Ferreira, L.~{Garrido-da-Silva}, and I.~S.
  Labouriau}, {\em Stability of cycles in a game of
  {R}ock-{S}cissors-{P}aper-{L}izard-{S}pock}, SIAM J. Appl. Dyn. Syst., 21
  (2022), pp.~2393--2431, \url{https://doi.org/10.1137/21M1435215}.

\bibitem{Castro2024}
{\sc S.~B. S.~D. Castro, A.~M.~J. Ferreira, and I.~S. Labouriau}, {\em
  Stability of cycles and survival in a jungle game with four species}, Dyn.
  Syst.,  (2024), pp.~1--19,
  \url{https://doi.org/10.1080/14689367.2024.2307515}.

\bibitem{Castro2025b}
{\sc S.~B. S.~D. Castro, C.~M. Postlethwaite, and A.~M. Rucklidge}, {\em
  Dataset for ``{V}isibility of heteroclinic networks''}.
\newblock University of Leeds Data Repository, 2025,
  \url{https://doi.org/10.5518/1652}.

\bibitem{Field1989}
{\sc M.~Field}, {\em Equivariant bifurcation theory and symmetry breaking}, J.
  Dyn. Diff. Equat., 1 (1898), pp.~369--421,
  \url{https://doi.org/10.1007/BF01048455}.

\bibitem{Glendinning1994}
{\sc P.~Glendinning}, {\em Stability, Instability and Chaos: An Introduction to
  the Theory of Nonlinear Differential Equations}, Cambridge University Press,
  Cambridge, 1994, \url{https://doi.org/10.1017/CBO9780511626296}.

\bibitem{GroothuizenDijkema2025}
{\sc D.~C. Groothuizen~Dijkema, V.~Kirk, C.~M. Postlethwaite, and A.~M.
  Rucklidge}, {\em Continuity of projected maps near heteroclinic networks in
  $\mathbb{R}^4$}, In preparation,  (2025).

\bibitem{Guckenheimer1988}
{\sc J.~Guckenheimer and P.~Holmes}, {\em Structurally stable heteroclinic
  cycles}, Math. Proc. Camb. Phil. Soc., 103 (1988), pp.~189--192,
  \url{https://doi.org/10.1017/S0305004100064732}.

\bibitem{Hofbauer1998}
{\sc J.~Hofbauer and K.~Sigmund}, {\em Evolutionary Games and Population
  Dynamics}, Cambridge University Press, Cambridge, 1998,
  \url{https://doi.org/10.1017/CBO9781139173179}.

\bibitem{Kirk2008}
{\sc V.~Kirk and A.~M. Rucklidge}, {\em The effect of symmetry breaking on the
  dynamics near a structurally stable heteroclinic cycle between equilibria and
  a periodic orbit}, Dyn. Syst. Int. J., 23 (2008), pp.~43--74,
  \url{https://doi.org/10.1080/14689360701709088}.

\bibitem{Kirk1994}
{\sc V.~Kirk and M.~Silber}, {\em A competition between heteroclinic cycles},
  Nonlinearity, 7 (1994), pp.~1605--1621,
  \url{https://doi.org/10.1088/0951-7715/7/6/005}.

\bibitem{Krupa1997}
{\sc M.~Krupa}, {\em Robust heteroclinic cycles}, J. Nonlin. Sci., 7 (1997),
  pp.~129--176, \url{https://doi.org/10.1007/BF02677976}.

\bibitem{Krupa1995}
{\sc M.~Krupa and I.~Melbourne}, {\em Asymptotic stability of heteroclinic
  cycles in systems with symmetry}, Ergod. Theory Dyn. Syst., 15 (1995),
  pp.~121--147, \url{https://doi.org/10.1017/S0143385700008270}.

\bibitem{Kuznetsov1998}
{\sc Y.~Kuznetsov}, {\em Elements of Applied Bifurcation Theory, Second
  Edition}, vol.~112 of Applied Mathematical Sciences, Springer-Verlag, New
  York, 1998, \url{https://doi.org/10.1007/978-3-031-22007-4}.

\bibitem{Leine2010}
{\sc R.~I. Leine}, {\em The historical development of classical stability
  concepts: {L}agrange, {P}oisson and {L}yapunov stability}, Nonlinear Dyn., 59
  (2010), pp.~173--182, \url{https://doi.org/10.1007/s11071-009-9530-z}.

\bibitem{Liapounoff1907}
{\sc A.~Liapounoff}, {\em Probl\`eme g\'en\'eral de la stabilit\'e du
  mouvement}, Ann. Fac. Sci. Toulouse (2), 9 (1907), pp.~203--474,
  \url{http://www.numdam.org/item/AFST_1907_2_9__203_0/}.

\bibitem{Mawhin1994}
{\sc J.~Mawhin}, {\em The centennial legacy of {P}oincar{\'e} and {L}yapunov in
  ordinary differential equaitons}, Rend. Circ. Mat. Palermo (2) Suppl., 34
  (1994), pp.~9--46,
  \url{https://www.researchgate.net/publication/242012997_The_centennial_legacy_of_Poincare_and_Lyapunov_in_ordinary_differential_equations}.

\bibitem{Melbourne1989a}
{\sc I.~Melbourne}, {\em Intermittency as a codimension-three phenomenon}, J.
  Dyn. Diff. Equat., 1 (1989), pp.~347--367,
  \url{https://doi.org/10.1007/BF01048454}.

\bibitem{Melbourne1991}
{\sc I.~Melbourne}, {\em An example of a non-asymptotically stable attractor},
  Nonlinearity, 4 (1991), pp.~835--844,
  \url{https://doi.org/10.1088/0951-7715/4/3/010}.

\bibitem{Milnor1985}
{\sc J.~Milnor}, {\em On the concept of attractor}, Commun. Math. Phys., 99
  (1985), pp.~177--195, \url{https://doi.org/10.1007/BF01212280}.

\bibitem{Podvigina2012}
{\sc O.~Podvigina}, {\em Stability and bifurcations of heteroclinic cycles of
  type ${Z}$}, Nonlinearity, 25 (2012), pp.~1887--1917,
  \url{https://doi.org/10.1088/0951-7715/25/6/1887}.

\bibitem{Podvigina2023}
{\sc O.~Podvigina}, {\em Behaviour of trajectories near a two-cycle
  heteroclinic network}, Dyn. Syst. Int. J., 38 (2023), pp.~576--596,
  \url{https://doi.org/10.1080/14689367.2023.2225463}.

\bibitem{Podvigina2020}
{\sc O.~Podvigina, S.~B. S.~D. Castro, and I.~S. Labouriau}, {\em Asymptotic
  stability of robust heteroclinic networks}, Nonlinearity, 33 (2020),
  pp.~1757--1788, \url{https://doi.org/10.1088/1361-6544/ab6817}.

\bibitem{Poincare1885}
{\sc H.~Poincar{\'e}}, {\em Sur les equations lin\'eaires aux diff\'erentielles
  ordinaires et aux diff\'erences finies}, Am. J. Math., 7 (1885), p.~203,
  \url{https://doi.org/10.2307/2369270}.

\bibitem{Postlethwaite2022}
{\sc C.~M. Postlethwaite and A.~M. Rucklidge}, {\em Stability of cycling
  behaviour near a heteroclinic network model of
  {R}ock--{P}aper--{S}cissors--{L}izard--{S}pock}, Nonlinearity, 35 (2022),
  p.~1702, \url{https://doi.org/10.1088/1361-6544/ac3560}.

\bibitem{Roque2011}
{\sc T.~Roque}, {\em Stability of trajectories from {P}oincar{\'e} to
  {B}irkhoff: approaching a qualitative definition}, Arch. Hist. Exact Sci., 65
  (2011), pp.~295--342, \url{https://doi.org/10.1007/s00407-011-0079-0}.

\end{thebibliography}


\end{document}